\documentclass[12pt]{article}
\usepackage[reqno]{amsmath}
\usepackage[all]{xy}
\usepackage[mathscr]{euscript}
\usepackage{amssymb}
\usepackage{amscd}
\usepackage{textcomp}
\usepackage{amsmath}
\usepackage{color}
\usepackage{a4wide}
\newtheorem{thm}{Theorem}

\newtheorem{lem}{Lemma}

\newtheorem{rem}{Remark}

\newcommand{{{\Ge}}}{\operatorname{G}}

\newcommand{\de}{\operatorname{d}}

\newcommand{\M}{{\cal M}}
\newcommand{\eeM}{{\cal M}^{\cdot \cdot}}
\newcommand{\eM}{{\cal M}^{\cdot}}

\newcommand{\B}{{\cal B}}

\newcommand{\F}{{\cal F}}

\newcommand{\ku}{\textit}
\newcommand{\ca}{\mathcal}

\title{A Characterization of the Poisson Process revisited}

\author{Benjamin Nehring\footnote{Fakult\"at f\"ur Mathematik der Ruhr-Universit\"at Bochum, Universit\"atsstra\ss{}e 150,
\newline 44780 Bochum, Germany, e-mail: {\sf Benjamin.Nehring@ruhr-uni-bochum.de}
}
}

\begin{document}

\maketitle

\begin{abstract}
We show that the splitting-characterization of the Poisson point process is an immediate consequence of the Mecke-formula. 
\end{abstract}


%%%%%%%%%%%%%%%%%%%%%%%%%%%%%%%%%%%%%%
\section{Notation}
%%%%%%%%%%%%%%%%%%%%%%%%%%%%%%%%%%%%%%%%%%%%

Let $X$ be a Polish space, $\B(X)$ resp. $\B_0(X)$ denote the Borel resp. bounded Borel sets. $\M(X)$ is the space of locally finite measures on X, i.e., Radon measures on $X$, which is Polish for the vague topology. $\eeM(X)$ denotes the closed and thereby measurable subspace of Radon point measures and $\eM(X)$ denotes the measurable subspace of simple Radon point measures. A law $P \in \ca{P}(\eeM(X))$ on $\eeM(X)$ resp. $\eM(X)$ is called point process resp. simple point process. The first moment measure of a point process $P$ will be denoted by 
\[
\nu_P(B) =  \int\limits_{\eeM(X)} P(\de \mu) \, \mu(B), \quad B\in \B(X).
\]
We say that a point process $P$ is of first order if $\nu_P \in \M(X)$. By $U$ we denote the set of non negative bounded measurable test functions on $X$ with a support contained in some bounded set. Remark that for $f\in U$, $\zeta_f: \mu \mapsto \mu(f)$ is a well defined measurable function on $\eeM(X)$ and let $\ca{L}_P(f) = P(e^{-\zeta_f})$, $f\in U$, be the Laplace transform of a point process $P$. 
Furthermore we let $T^\mu_q$ be the independent q-thinning $q\in (0,1)$ of a Radon point measure $\mu \in \eeM(X)$, that is 
\[
T^\mu_q = \underset{x\in \mu}{\ast} ((1-q)\, \delta_{\mathbf{0}} + q \, \delta_{\delta_x})^{\ast \mu(\{x\})}.
\]
Here $\ast$ denotes ordinary convolution of point processes and $\mathbf{0}$ is the zero measure on $X$. "$x\in \mu$" means that the convolution is taken over the set $\{y\in X \, | \, \mu(\{y\}) > 0 \}$. The independent q-thinning $\Gamma_q(P)$ of a point process $P$ can now be defined as 
\[
\Gamma_q(P) = \int\limits_{\eeM(X)} P(\de \mu) \,\, T^\mu_q.
\]

%%%%%%%%%%%%%%%%%%%%%%%%%%%%%%%%%%%%%%%%%%%%%%
\section{Introduction}
%%%%%%%%%%%%%%%%%%%%%%%%%%%%%%%%%%%%%%%%%%%%%%%%%

Having observed a realization $\nu \in \eeM(X)$ of an independent q-thinning of a point process $P$ one can ask the following question: What is the distribution of deleted point configurations given the realization $\nu \in \eeM(X)$? This conditional probability will be called splitting kernel and will be denoted by $\Upsilon^\nu_q(P)$ in the sequel (see also section 6.3 in \cite{N2}). In case $P$ is a finite point process, that is $P$ is concentrated on the set of finite Radon point measures, Karr obtained in \cite{Karr} a representation of the splitting kernel in terms of the reduced Palm distributions $P^!_{\delta_{x_1}+\ldots+\delta_{x_n}}$, $x_1,\ldots,x_n \in X$, $n\geq 1$, of $P$. That is 
\[
\Upsilon^\nu_q(P)(\varphi) = \frac{1}{\int (1-q)^{\mu(X)} P^!_\nu(\de \mu)} \int \varphi(\mu)\,(1-q)^{\mu(X)} P^!_\nu(\de \mu),
\]
where $\varphi\in F_+$ is some non negative measurable test function on the space of finite point measures and $\nu$ denotes a finite point measure (see also proposition 6.3.5 in \cite{N2}). To obtain a representation for $\Upsilon^\nu_q(P)$ in case $P$ is a general point process seems to be an open problem. \\

In this note we want to prove the following: Assume you have a point process $P$ such that $\Upsilon^\nu_q(P)$ does not depend on the observed point configuration $\nu \in \eeM(X)$. That is, there is some point process $Q_q$ such that $\Upsilon^\nu_q(P) = Q_q$ for all $\nu \in \eeM(X)$. Then $P$ can only be a Poisson point process. This result is a corollary of Fichtner's main theorem (Satz 1) in \cite{Fichtner}. Fichtner's arguments were quite involved so Assun\c{c}\~{a}o and Ferrari \cite{Assu}  gave a simpler proof of the result (in the present setting) using a characterization of the Poisson distribution and the fact that a simple point process is determined by its avoidance function. Note that in \cite{Assu} the result is stated for general point processes (meaning elements of $\ca{P}(\eeM(X))$) but i.e. Brown and Xia have shown in \cite{Brown} that one can in general not conclude that the point process is Poisson if its counting variables $\zeta_B$, $B\in \B_0(X)$, are Poisson distributed. In fact prior to Assun\c{c}\~{a}o and Ferrari, Shengwu and Jiagang have stated the result more precisely in \cite{SJ} for simple point processes by using that a simple point process is determined by the distribution of its counting variables. Since in the present setting we do not want to exclude the possibility of multiple points we have to resort to different techniques. The most important one will be Mecke's characterization of the Poisson point process (Satz 3.1 in \cite{M1}). Moreover in section 3.1 it is shown that the same techniques also cover the case of "multi-splitting".\\
Let us introduce the notion of a Papangelou kernel  (also sometimes called conditional intensity) $\pi$ of a point process $P$. For a detailed description we refer to \cite{Z}, where also the historic development is explained. $\pi$ is a kernel from $\eeM(X)$ to $\M(X)$, that is for any $\mu \in \eeM(X)$ we have $\pi(\mu,\de x) \in \M(X)$, so that $P$ satisfies the equation
\[
C_P(h) := \int\limits_{\eeM(X)} \int\limits_X h(x,\mu)\, \mu(\de x) P(\de \mu) = \int\limits_{\eeM(X)} \int\limits_X h(x,\mu+\delta_x) \, \pi(\mu,\de x) P(\de \mu),
\] 
for all non negative measurable test functions $h$ on the product space $X \times \eeM(X)$. In the first equation the definition of the Campbell measure of $P$ is provided.\\
One direction of Mecke's result can now be formulated as follows: Assume $P$ is a point process whose Papangelou kernel $\pi(\mu,\de x)$ does not depend on $\mu \in \eeM(X)$, that is, there is some $\varrho \in \M(X)$ such that $\pi(\mu,\cdot) = \varrho$ for all $\mu \in \eeM(X)$ then $P$ is a Poisson point process with first moment measure given by $\varrho$. \\
\\
For the proof of the characterization result we will need as a main lemma:
\begin{lem}\label{rem1}
For all $q\in (0,1)$ and $\mu \in \eeM(X)$, $T^\mu_q$ has a Papangelou kernel given by
\[
\pi(\kappa,\de x) = \frac{q}{1-q} \,(\mu - \kappa)(\de x), \quad \kappa \in \eeM(X).
\]
Note that $T^\mu_q$ realizes only sub configurations of $\mu$ so $\mu - \kappa$ is $T^\mu_q$ - a.s. $[\kappa]$ in $\eeM(X)$.
\end{lem}
A proof can be found in the scholion "The P\'olya difference process".\\
Furthermore we need the so called splitting law $S_q(P)$ of a point process $P$
\[
S_q(P)(h) = \int\limits_{\eeM(X)}\int\limits_{\eeM(X)} P(\de \mu) T^\mu_q(\de \nu) \, h(\nu,\mu-\nu),
\]
where $h$ is some non negative measurable test function on $\eeM(X)\times \eeM(X)$. So $S_q(P)$ is a law on $\eeM(X)\times \eeM(X)$, which realizes tuples $(\nu,\eta)$ such that $\nu$ is the point configuration which survived the thinning and $\eta$ is the collection of deleted points. The marginal laws of $S_q(P)$ are given by
\[
S_q(P)(\varphi \otimes \mathbf{1}) = \Gamma_q(P)(\varphi) \text{ and } S_q(P)(\mathbf{1}\otimes \varphi) = \Gamma_{1-q} (P)(\varphi),
\]
where $\varphi$ is a non negative test function on $\eeM(X)$ and $\mathbf{1}$ denotes the function, which is constantly one. Thus for any $N\in\B(\eeM(X))$ we have that $S_q(P)(\cdot \times N)$ is absolutely continuous to $\Gamma_q(P)$. Therefore by the theory of disintegration we obtain the existence of the splitting kernel $\Upsilon^\nu_q(P)$, that is
\[
S_q(P)(\de \nu  \de \eta) = \Gamma_q(P)(\de \nu)\,  \Upsilon^\nu_q(P)(\de \eta).
\]
%%%%%%%%%%%%%%%%%%%%%%%%%%%%%%%%%%%%%%%%%%%%%%%%%%%%%%%%%%%%%%%%%%%%%%%%%%%%%%
\section{A Characterization }
%%%%%%%%%%%%%%%%%%%%%%%%%%%%%%%%%%%%%%%%%%%%%%%%%%%%%%%%%%%%%%%%%%%%%%%%%%%%%%

We are now ready to state the result.

\begin{thm}
Let $P$ be a point process of first order. Then $P$ is a Poisson point process if and only if the splitting law factorizes into its marginals, that is
\[
(\F) \quad \quad S_q(P) = \Gamma_q(P) \otimes \Gamma_{1-q}(P),
\]
for some $q\in(0,1)$.
\end{thm}
\ku{Proof}.
Let us denote by $\Pi_\varrho$ the Poisson point process with first moment measure $\varrho \in \M(X)$. First assume $P = \Pi_\varrho$. To establish the identity $(\F)$ it suffices to show the equality of both sides on the class of test functions $h = e^{-\zeta_f}\otimes e^{-\zeta_g}$, where $f,g \in U$.  We have 
\[
S_q(\Pi_\varrho)(h) = \int \Pi_\varrho(\de \mu)\, e^{-\mu(g)}\, T^\mu_q(e^{\zeta_{g-f}}).
\]
Note that by definition of $T^\mu_q$
\[
T^\mu_q(e^{\zeta_{g-f}}) = \prod_{x\in \mu} \left(1-q + q e^{(g-f)(x)}\right)^{\mu(\{x\})} =\exp(\mu(\log(1-q + q \,e^{g-f})))
\]
and remark that only finitely many factors in the above product are different from one. So we obtain $S_q(\Pi_\varrho)(h) = \ca{L}_{\Pi_\varrho}(v)$, where $v=g - \log(1-q+q\, e^{g-f})$. One straightforwardly checks that $v\in U$. In fact $v = -\log((1-q)\, e^{-g} + q\, e^{-f} )$. By using the representation of the Laplace transform of $\Pi_\varrho$ one obtains
\[
 \ca{L}_{\Pi_\varrho}(v) =  \ca{L}_{\Pi_{q \varrho}}(f)\, \ca{L}_{\Pi_{(1-q) \varrho}}(g),
\] 
which establishes the identity $(\F)$, since it is well known that $\Gamma_s(\Pi_\varrho) = \Pi_{s \varrho}$ for any $s\in (0,1)$.\\
Assume now from the contrary that $P$ solves $(\F)$ for some $q\in (0,1)$. Let us compute the Campbell measure of $\Gamma_q(P)$. Take $h$ to be a non negative measurable test function on $X\times \eeM(X)$ then we have 
\begin{align*}
&C_{\Gamma_q(P)}(h)  \overset{(i)}{=} \int P(\de \mu) C_{T^\mu_q}(h) 
 \overset{(ii)}{=} \int P(\de \mu) T^\mu_q(\de \kappa) \, \frac{q}{1-q} (\mu-\kappa)(\de x) \, h(x,\kappa+\delta_x)\\
& \overset{(iii)}{=} \int S_q(P)(\de \kappa \de \eta) \, \frac{q}{1-q} \eta(\de x) \, h(x,\kappa+\delta_x)\\
& \overset{(iv)}{=} \int \Gamma_q(P)(\de \kappa) \Gamma_{1-q}(P)(\de \eta) \, \frac{q}{1-q} \eta(\de x) \, h(x,\kappa+\delta_x)\\
& \overset{(v)}{=} \int \Gamma_q(P)(\de \kappa) \frac{q}{1-q} \nu_{\Gamma_{1-q}(P)}(\de x)  \, h(x,\kappa+\delta_x)
 \overset{(vi)}{=} \int \Gamma_q(P)(\de \kappa)\, q\, \nu_P(\de x) \, h(x,\kappa+\delta_x).
\end{align*}
$(i)$ follows by definition of $\Gamma_q(P)$. $(ii)$ is due to lemma \ref{rem1} in the introductory section. $(iii)$ follows by definition of the splitting law $S_q(P)$. Since $P$ is assumed to satisfy $(\F)$ $(iv)$ holds. In $(v)$ the definition of the first moment measure has been used. Finally $(vi)$ holds true because $\nu_{\Gamma_{1-q}(P)} = (1-q)\,  \nu_P$.
So by Mecke's characterization (Satz 3.1 in \cite{M1}) it follows that $\Gamma_q(P) = \Pi_{q\, \nu_P}$. Lemma 9 in \cite{M2} states that $\Gamma_q : \ca{P}(\eeM(X)) \rightarrow \ca{P}(\eeM(X))$ is an injective mapping (See also lemma 2 in \cite{SJ}. Their argument is also valid for arbitrary point processes (not only for the simple ones)). Therefore $P = \Pi_{\nu_P}$.\\
$\square$

%%%%%%%%%%%%%%%%%%%%%%%%%%%%%%%%%%%%%%%%%%%%%%%%%%%%%%%%%%%%%%%%%%%%%%%%%%%%%%%%%%%%%%%%%%%%%%%%%%%%%%%%%%%%%%%%%%%%%%%%%%%%%%%%%%%%%%%%%%%%%%%%%%%%%%%%%%%%%%%%%%%%%%%%%%%%%%%%%%%%%%%%
\subsection{A Generalization}
%%%%%%%%%%%%%%%%%%%%%%%%%%%%%%%%%%%%%%%%%%%%%%%%%%%%%%%%%%%%%%%%%%%%%%%%%%%%%%%%%%%%%%%%%%%%%%%%%%%%%%%%%%%%%%%%%%%%%%%%%%%%%%%%%%%%%%%%%%%%%%%%%%%%%%%%%%%%%%%%%%%%%%%%%%%%%%%%%%%%%%%%

As it was done in \cite{Fichtner} one can not only obtain a characterization of the Poisson point process if one partitions the randomly realized point measure $\mu \in \eeM(X)$ in two sub configurations but also in $n\geq 2$. Let us investigate the following mechanism: Let $\mu \in \eeM(X)$ and $s_1,\ldots,s_{n-1} \in (0,1)$. Then we define
 \begin{align*}
T^\mu_{(s_1,\ldots,s_{n-1})}(h) = \int T^\mu_{s_1}(\de \kappa_1)  T^{\mu-\kappa_1}_{s_2}(\de \kappa_2) & \ldots  T^{\mu-(\kappa_1+\ldots+\kappa_{n-2})}_{s_{n-1}}(\de \kappa_{n-1})\\
&    h(\kappa_1,\kappa_2,\ldots, \kappa_{n-1},\mu-(\kappa_1+\ldots+\kappa_{n-1})),
\end{align*}
where $h\in F_+$ is some non negative measurable test function on the space $\overset{n}{\underset{j=1}{\times}} \eeM(X)$. So $T^\mu_{(s_1,\ldots,s_{n-1})}$ partitions the point configuration $\mu$ in $n$ parts 
\[
\mu = \kappa_1 +\ldots + \kappa_n \text{ with } \kappa_m \in \eeM(X), \quad 1\leq m \leq n.
\]
The probability that a point of $\mu$ belongs to the $m$-th collection of points is given by
\[
(1-s_1) (1-s_2)\ldots (1-s_{m-1}) s_m, \quad 1\leq m \leq n.
\]
We are setting here $s_n = 1$. Let now $q_1, \ldots, q_n \in (0,1)$, such that they sum up to one. If we now want that a point of $\mu$ has probability $q_m$ to belong to the $m$-th collection of points we have to adjust the $s_m$'s appropriately. In fact the right choice is
\begin{equation}\label{eq1}
s_m = \frac{q_m}{q_m+\ldots +q_{n}}, \quad 1\leq m \leq n.
\end{equation}
We can now introduce the multi-splitting law
\[
S^n_{q} (P) = \int\limits_{\eeM(X)} P(\de \mu)\, T^\mu_{(s_1,\ldots,s_{n-1})},
\]
where the $s_m$'s have been chosen as in (\ref{eq1}) and $q$ denotes the law on $\{1,\ldots,n\}$ such that $q(\{m\})=q_m$. We are now ready to state the result.

\bigskip
%\bigskip

\begin{thm}
Let $P$ be a point process of first order and $q_1, \ldots, q_n \in (0,1)$ such that 
\[
q_1+\ldots+q_n=1.
\]
Then $P$ is a Poisson point process if and only if 
\[
(\F_n) \quad \quad S^n_{q} (P) = \Gamma_{q_1}(P) \otimes \ldots \otimes \Gamma_{q_n}(P).
\] 
\end{thm}

\medskip

\begin{rem}
Note that the first orderness of $P$ is necessary to apply Mecke's characterization.
\end{rem}

%%%%%%%%%%%%%%%%%%%%%%%%%%%%%%%%%% Neue Seite %%%%%%%%%%%%%%%%%%%%%%%%%%%%%%%%%%%%%%%%

\newpage

\ku{Proof.} We will show by induction that the Poisson point process $\Pi_\varrho$, $\varrho \in \eeM(X)$ satisfies $(\F_n)$ for all $n\geq 2$ and distributions $q$. We will use that the multi-splitting law $S^n_q(P)$ is given recursively: Let $w$ be the distribution on $\{1,\ldots,n-1\}$ such that $w_j = q_j$, $1\leq j \leq n-2$ and $w_{n-1} = q_{n-1} + q_n$. Then for any $h\in F_+$ and point process $P$
\begin{align*}
S^n_{q} (P)(h) = \int S^{n-1}_w(P) (\de \kappa_1,\ldots, \de \kappa_{n-1}) & \,  T^{\kappa_{n-1}}_{\frac{q_{n-1}}{q_{n-1}+q_n}}(\de \eta) \,\,
 h(\kappa_1,\ldots,\kappa_{n-2},\eta, \kappa_{n-1} -\eta).
\end{align*}
So for the Poisson point process $\Pi_\varrho$ we have by induction
\begin{align*}
S^n_{\sigma} (\Pi_\varrho)(h) & =  \int  \Gamma_{q_1}(\Pi_\varrho)(\de \kappa_1)  \ldots  \Gamma_{q_{n-2}}(\Pi_\varrho)(\de \kappa_{n-2}) \Gamma_{q_{n-1}+q_{n}}(\Pi_\varrho)(\de \kappa_{n-1})\\
&   \quad \quad \quad 	 \quad \quad \quad \quad		 T^{\kappa_{n-1}}_{\frac{q_{n-1}}{q_{n-1}+q_n}}(\de \eta) \,  h(\kappa_1,\ldots,\kappa_{n-2},\eta, \kappa_{n-1} -\eta)\\
 &= \int  \Pi_{q_1\varrho}(\de \kappa_1)  \ldots  \Pi_{q_{n-2}\varrho}(\de \kappa_{n-2}) \\
&  \quad \quad 	 \quad 		S_{\frac{q_{n-1}}{q_{n-1}+q_n}}(\Pi_{(q_{n-1}+q_{n}) \varrho})(\de \eta, \de \zeta) \, h(\kappa_1,\ldots,\kappa_{n-2},\eta, \zeta)\\
&= \Pi_{q_1\varrho}\otimes \ldots \otimes \Pi_{q_n \varrho} (h)
\end{align*}
The last equality holds since $(\F_2)$ has already been established for the Poisson processes in theorem 1.\\
Assume now from the contrary that $P$ is a solution to $(\F_n)$. We will follow the proof of theorem 1 and compute the Campbell measure of $\Gamma_{q_1}(P)$
\begin{align*}
&C_{\Gamma_{q_1}(P)}(h)  \overset{(i)}{=} \int P(\de \mu) T^\mu_{q_1}(\de \kappa_1) \, \frac{q_1}{1-q_1} (\mu-\kappa_1)(\de x) \, h(x,\kappa_1+\delta_x)\\
& \overset{(ii)}{=} \int P(\de \mu)  T^\mu_{s_1}(\de \kappa_1)  T^{\mu-\kappa_1}_{s_2}(\de \kappa_2) \ldots  T^{\mu-(\kappa_1+\ldots+\kappa_{n-2})}_{s_{n-1}}(\de \kappa_{n-1}) \\
&     			\quad\quad\quad			 \frac{q_1}{1-q_1} (\kappa_2+\ldots+\kappa_{n-1}+ (\mu - (\kappa_1+\ldots+\kappa_{n-1}))(\de x) \, h(x,\kappa_1+\delta_x)\\
& \overset{(iii)}{=} \int S^n_q(P)(\de \kappa_1\ldots \de \kappa_n) \, \, \frac{q_1}{1-q_1} (\kappa_2+\ldots+\kappa_{n})(\de x) \, h(x,\kappa_1+\delta_x)   \\
& \overset{(iv)}{=} \int \Gamma_{q_1}(P)(\de \kappa_1)\ldots \Gamma_{q_n}(P)(\de \kappa_n)  \, \, \frac{q_1}{1-q_1} (\kappa_2+\ldots+\kappa_{n})(\de x) \, h(x,\kappa_1+\delta_x) \\
& \overset{(v)}{=} \int \Gamma_{q_1}(P)(\de \kappa_1)  \frac{q_1}{1-q_1} (\nu_{\Gamma_{q_2}(P)}+\ldots+\nu_{\Gamma_{q_n}(P)})(\de x) \, h(x,\kappa_1+\delta_x).\\
& \overset{(vi)}{=} \int \Gamma_{q_1}(P)(\de \kappa_1) \,   q_1 \nu_P(\de x) \, h(x,\kappa_1+\delta_x).
\end{align*}
$(i)$ follows exactly as in the proof of theorem 1. In $(ii)$ an appropriate zero measure has been added and note that $q_1 = s_1$. $(iii)$ is due to the definition of $S^n_q(P)$. In $(iv)$ $(\F_n)$ has been used. The equality $(v)$ holds by definition of the first moment measure. Finally $(vi)$ holds since, as remarked earlier, $\nu_{\Gamma_{q_j}(P)} = q_j \nu_P$. So we conclude as in the proof of theorem 1. $\square$

%%%%%%%%%%%%%%%%%%%%%%%%%%%%%%%%%%%%%%%%%%%%%%%%%%%%%%%%%%%%%%%%%%%%%%%%%%%%%%%%%%%%%%%%%%%%%%%%%%%%%%%%%%%%%%%%%%%%%%%%%%%%%%%%%%%%%%%%%%%%%%%%%%%%%%%%%%%%%%%%%%%%%%%%%%%%%%%%%%%%%%%%
\section{Scholion: The P\'olya difference process}
%%%%%%%%%%%%%%%%%%%%%%%%%%%%%%%%%%%%%%%%%%%%%%%%%%%%%%%%%%%%%%%%%%%%%%%%%%%%%%%%%%%%%%%%%%%%%%%%%%%%%%%%%%%%%%%%%%%%%%%%%%%%%%%%%%%%%%%%%%%%%%%%%%%%%%%%%%%%%%%%%%%%%%%%%%%%%%%%%%%%%%%%
\ku{
In this section we will give a proof of lemma \ref{rem1}. The section is named "P\'olya difference process", because the P\'olya difference process $P^-_{z,\mu}$, $z\in (0,\infty)$, $\mu\in \eeM(X)$ as introduced in \cite{NZ} was shown to have a Papangelou kernel given by 
\[
\pi(\kappa,\de x) = z\, (\mu-\kappa)(\de x), \quad \kappa \in \eeM(X).
\]  
Furthermore we have $T^\mu_q = P^-_{\frac{q}{1-q},\mu}$. So the result of lemma \ref{rem1} is clear. However we will give here an independent proof.\\ \\
}
\ku{Proof of lemma \ref{rem1}}. 
Let $a\in X$. Then it is immediately checked that 
\[
\left(\Sigma' \right) \quad C_{T^{\delta_a}_q}(h) = \iint h(x, \eta + \delta_x) \, \frac{q}{1-q} (\delta_a - \eta)(\de x) T^{\delta_a}_q(\de \eta)
\]
holds for any non negative measurable test function $h$ on $X \times \eeM(X)$, since $T^{\delta_a}_q$ realizes either the zero measure with probability $1-q$ or a point $\delta_a$ with probability $q$. Let now $\mu \in \eeM(X)$ be arbitrary. 
\begin{align*}
& C_{T^{\mu}_q}(h)  \overset{(i)}{=} \int h\bigg(y, \sum\limits_{x\in \mu} \sum\limits_{j=1}^{\mu(\{x\})} \eta^j_x\bigg) \sum\limits_{a\in \mu} \sum\limits_{i=1}^{\mu(\{a\})} \eta^i_a (\de y)
\prod\limits_{x\in \mu} \prod\limits_{j=1}^{\mu(\{x\})} T^{\delta_x}_q(\de \eta^j_x) \\
& \overset{(ii)}{=} \sum\limits_{a\in \mu} \sum\limits_{i=1}^{\mu(\{a\})}  \int h\bigg(y, \sum\limits_{x\in \mu} \sum\limits_{j=1}^{\mu(\{x\})} \eta^j_x\bigg) \eta^{i}_a(\de y) \prod\limits_{x\in \mu} \prod\limits_{j=1}^{\mu(\{x\})} T^{\delta_x}_q(\de \eta^j_x)\\
& \overset{(iii)}{=} \sum\limits_{a\in \mu} \sum\limits_{i=1}^{\mu(\{a\})}  \int h\bigg(y, \sum\limits_{x\in \mu} \sum\limits_{j=1}^{\mu(\{x\})} \eta^j_x +\delta_y\bigg) \, \frac{q}{1-q} (\delta_a - \eta^{i}_a)(\de y)
\prod\limits_{x\in \mu} \prod\limits_{j=1}^{\mu(\{x\})} T^{\delta_x}_q(\de \eta^j_x)\\
& \overset{(iv)}{=} \int h\bigg(y, \sum\limits_{x\in \mu} \sum\limits_{j=1}^{\mu(\{x\})} \eta^j_x +\delta_y\bigg)  \frac{q}{1-q} (\mu - \sum\limits_{a\in \mu} \sum\limits_{i=1}^{\mu(\{a\})} \eta^{i}_a)(\de y)
\prod\limits_{x\in \mu} \prod\limits_{j=1}^{\mu(\{x\})} T^{\delta_x}_q(\de \eta^j_x)\\
& \overset{(v)}{=} \int h(y, \eta +\delta_y)\, \frac{q}{1-q} (\mu - \eta)(\de y)\, T^\mu_q(\de \eta).
\end{align*}
In $(i)$ the definition of the Campbell measure and the fact that 
\[
T^\mu_q = \underset{x\in \mu}{\ast} \left( T^{\delta_x}_q \right) ^{\ast \mu(\{x\})}
\]
has been used. In $(ii)$ summation and integration have been interchanged. $(iii)$ is due to $\left(\Sigma' \right)$. Steps $(iv)$ and $(v)$ are steps $(i)$ and $(ii)$ in the reversed order. $\square$

\newpage

\textbf{\ku{\large{Acknowledgement}}}
\newline
\newline
I would like to thank the Sonderforschungsbereich - TR 12 of the Deutsche\\ Forschungsgemeinschaft for financial support.

 %%%%%%%%%%%%%%%%%%%%%%%%%%%%%

%%%%%%%%%%%%%%%%%%%%%%%%%%%%%%%%%%

\end{document}